\documentclass{amsart}

\usepackage{amssymb}
\usepackage{graphicx}
\usepackage{flafter}

\newtheorem{theorem}{Theorem}[section]

\newtheorem{lemma}[theorem]{Lemma}

\theoremstyle{definition}

\theoremstyle{remark}

\numberwithin{equation}{section}

\newcommand\blfootnote[1]{%
  \begingroup
  \renewcommand\thefootnote{}\footnote{#1}%
  \addtocounter{footnote}{-1}%
  \endgroup
}

\begin{document}

\title[Area compression of homeomorphisms with integrable distortion]{Stretching multifractal spectra and area compression of homeomorphisms with integrable distortion in higher dimensions}
\keywords{Mappings of finite distortion, rotation, integrable distortion. \newline \indent The author was financially supported by the V\"ais\"al\"a Foundation and by The Centre of Excellence in Analysis and Dynamics Research (Academy of Finland, decision 271983)}

\author{Lauri Hitruhin}
\address{University of Helsinki, Department of mathematics and Statistics, P-O. Box 68, FIN-00014 University of Helsinki, Finland}
\email{lauri.hitruhin@helsinki.fi}


\begin{abstract}
We consider homeomorphisms with integrable distortion in higher dimensions and sharpen the previous bound for area compression, which was presented by Clop and Herron in \cite{CH}. Our method relies on developing sharp bounds for the stretching multifractal spectra of these mappings.
\end{abstract}

 \maketitle

\section{Introduction}\label{INTRO}
Pointwise stretching of homeomorphisms with integrable distortion has been studied by Koskela and Takkinen in the planar case and by Clop and Herron in higher dimensions, see \cite{KT} and \cite{CH}. They proved that given an arbitrary  homeomorphic mapping  $f: \mathbb{R}^n \to \mathbb{R}^n $ with $p$-integrable distortion, where $p>n-1$, the pointwise stretching satisfies 
\begin{equation}\label{CHyht}
|f^{-1}(z)-f^{-1}(x)| \leq c_{p,n,K_f,G} \left(\frac{1}{\log\left( \frac{1}{|z-x|} \right)} \right)^{\frac{p(n-1)}{n}},
\end{equation}
where $|z-x|<1$ and $G=f^{-1}(B(z,1))$. Moreover, they verified that the exponent $\frac{p(n-1)}{n}$ in \eqref{CHyht} is optimal. \newline \newline
The sharp pointwise bound \eqref{CHyht} provides a starting point for the study of the stretching multifractal spectra, which measures the maximal size of a set in which these mappings can attain some predefined stretching. In the planar case this was done in \cite{H4}, and one of the main goals of this paper is to generalize this result to higher dimensions.
\begin{theorem}\label{MSPEKTRI}
Let $f: \mathbb{R}^n \to  \mathbb{R}^n$ be a homeomorphism with $p$-integrable distortion, where $p>n-1$, and fix $s\in (0,n)$. Furthermore, let  $A \subset   \mathbb{R}^n$ be the set of points for which there exists a sequence of  numbers  $|\lambda_{z,n}|\to 0$ such that
\begin{equation}\label{MSYHT}
\left|f^{-1}(z+\lambda_{z,n})-f^{-1}(z)\right| \geq c \left( \log\left( \frac{1}{|\lambda_{z,n}|} \right) \right)^{\frac{-p(n-1)}{n-s}},
\end{equation}
where $c>0$ is some fixed constant. Then the set $A$ satisfies 
\begin{equation*}
H^{h_a}(A)=0
\end{equation*}
for any $a> \frac{ps(n-1)}{n-s}$, where the gauge-function $h_a$ is defined by 
\begin{equation*}
h_a(r)=\left( \frac{1}{\log\left( \frac{1}{r} \right)} \right)^{a}.
\end{equation*}
\end{theorem}
See section \ref{PRE} for details regarding gauge-functions and generalized Hausdorff measures. \newline \newline 
Good understanding of the multifractal spectra of these mappings paves the way for the study of area compression. For quasiconformal mappings the sharp dimensional bounds for area compression were given by Astala in \cite{AS}. Further progress  in quasiconformal case, mostly dealing with replacing the Hausdorff dimension with measure, has been made by, for example, Astala, Clop, Mateu, Orobitg and Uriarte-Tuero, see \cite{ACM} and \cite{U}. While the quasiconformal case has been studied quite intensively one can also ask if similar bounds could be found for a more general family of mappings of finite distortion. \newline \newline 
 To this end, Clop and Herron in their article  \cite{CH} used the pointwise stretching bound \eqref{CHyht}  to estimate compression of small balls under $p$-integrable homeomorphisms. With this method they proved  that if  $f: \mathbb{R}^n \to \mathbb{R}^n$ is a homeomorphism with $p$-integrable distortion and $A \subset \mathbb{R}^n$ satisfies $H^{s}(A)>0$, then the image satisfies $H^{h}(f(A))>0$, where 
\begin{equation*}
h(r)=\left( \frac{1}{\log\left( \frac{1}{r} \right)} \right)^{\frac{ps(n-1)}{n}}.
\end{equation*}
Moreover, they constructed examples of homeomorphisms $f: \mathbb{R}^n \to \mathbb{R}^n$ with $p$-integrable distortion that can map a set $A \subset \mathbb{R}^n$, with $H^s(A)>0$, to a set which satisfies $H^{\bar{h}}(f(A))=0$, where 
\begin{equation*}
\bar{h}(r)=\left( \frac{1}{\log\left( \frac{1}{r} \right)} \right)^{\frac{ps(n-1)}{n-s}}.
\end{equation*}
As there was a gap left between the gauge functions $h$ and $\bar{h}$ Clop and Herron asked if  the result on area compression could be improved? \newline \newline
In the planar case we managed to do this, see \cite{H4}, using the stretching multifractal spectra, instead of the pointwise bound \eqref{CHyht}, to estimate compression of small balls. In this article our aim is to generalize this approach to higher dimensions.
\begin{theorem}\label{AREAC}
Let $s\in (0,n)$, $p>n-1$ and $f: \mathbb{R}^n \to \mathbb{R}^n$ be a homeomorphism with $p$-integrable distortion. Assume furthermore that 
\begin{equation*}
H^{s}(A)>0
\end{equation*}
for some set $A \subset \mathbb{R}^n$. Then 
\begin{equation*}
H^{h_a}(f(A))>0
\end{equation*}
whenever $a<\frac{ps(n-1)}{n-s}$.
\end{theorem}
Note, that Theorem \ref{AREAC} together with the examples constructed by Clop and Herron ensure that the gauge 
\begin{equation*}
h(r)=\left( \frac{1}{\log\left( \frac{1}{r} \right)} \right)^{\frac{ps(n-1)}{n-s}}
\end{equation*}
is indeed the critical one when measuring area compression.

\section{Prerequiseties}\label{PRE}
Let $\Omega \subset \mathbb{R}^n $ be a domain. We say that a homeomorphism $f: \Omega \to \mathbb{R}^n$ has finite distortion if the following conditions hold:
\begin{itemize}
\item $f\in W_{\text{loc}}^{1,1}(\Omega)$
\item $J_f(z)\in L^{1}_{\text{loc}}(\Omega)$
\item $|Df(z)|^n\leq J_f(z)K(z) \qquad \text{almost everywhere in $\Omega$},$
\end{itemize}
for a measurable function $K(z)\geq 1$, which is finite almost everywhere. The smallest such function is denoted by $K_f(z)$ and called the distortion of $f$. Here $Df(z)$ denotes the differential matrix of $f$ at the point $z$ and  $|Df(z)|$ is its operator norm,
whereas $J_f(z)$ is the Jacobian of the mapping $f$ at the point $z$. \newline \newline 
Such a mapping is said to have a $p$-integrable distortion, where $p\geq 1$, if 
\begin{equation*}
K_{f}(z)\in L_{\text{loc}}^{p}(\Omega).
\end{equation*}
For a detailed exposition of mappings of finite distortion see, for example, \cite{AIM} or \cite{HK}. \newline \newline
Our proof of Theorem \ref{MSPEKTRI} relies on estimates for the capacity of condensers. We briefly present here the results necessary for this paper, for a closer look on the topic we recommend, for example, \cite{VU}. \newline \newline 
Let $A \subset \mathbb{R}^n$ be open and $E\subset A$ compact, and call the pair $(E,A)$ a condenser. The $p$-capacity of a condenser is defined by
\begin{equation*}
\text{cap}_{p}(E,A)= \inf _{u} \int_{A} |\nabla u|^p \; dz,
\end{equation*}
where the infimum is taken over all continuous sobolev $W^{1,1}-$regular mappings with compact support in the set $A$ that satisfy $u(z)\geq 1$ when $z\in E$. Furthermore, standard approximation estimates let us assume that $u$ is a $C^{\infty}$ mapping that has compact support in the set $A$ and satisfies $0\leq u(z) \leq 1$ for all $z\in A$, see, for example, \cite{HKM}. We call these mappings admissible for the condenser $(E,A)$.
In our situation the open set $A$ will consist of finite number of disjoint  bounded domains. \newline \newline 
Capacity of a given set is usually impossible to calculate exactly, but for some trivial cases it is well known. For example, given $n\geq 2$ and $0<r<R<\infty$ we can calculate 
\begin{equation}\label{KAPTULOS}
\text{cap}_{n}\left( \overline{B(z,r)}, B(z, R) \right)=c_n \log^{1-n}\left( \frac{R}{r} \right),
\end{equation}
see, for example, \cite{VU}. In a more general setting we can estimate the capacity in the following way.
\begin{lemma}\label{KAPLEM}
Let $\Omega \subset \mathbb{R}^n$ be  a bounded domain and $E  \subsetneq \Omega$ a continuum. Then for any constant  $q$ such that $n-1<q<n$, we can estimate 
\begin{equation}\label{KAPYHT}
\text{cap}_{q}(E,\Omega)\geq c_{q,n}\left( \textrm{diam}(E) \right)^{n-q}.
\end{equation}
For the proof see \cite{KT} by Koskela and Takkinen in the planar case and \cite{CH} by Clop and Herron in the case $n>2$. 
\end{lemma}
When describing area compression of homeomorphisms with integrable distortion we need  more delicate scales than the classical Hausdorff measures. Instead we have to use more general Hausdorff gauge-functions, which are non-decreasing functions $h:(0,\infty) \to (0, \infty)$ such that $\lim_{r\to 0}h(r)=0$. We define the Hausdorff measure $H^h$ of a set $A\subset \mathbb{R}^n$ for each of these gauge-functions by 
\begin{equation*}
H^h(A)= \lim_{r \to 0}\left[ \inf \left\{ \sum_{i} h(\text{diam}(A_i)) \; : \; A\subset \bigcup_{i} A_i  \; \; \text{and} \; \; \text{diam}(A_i)<r \right\} \right].
\end{equation*}
It is well known that given an arbitrary set  $A \subset \mathbb{R}^n$ and any gauge-functions $h$ and $g$ we have the inequality
\begin{equation} \label{QEY}
H^{h}(A) \leq \limsup_{r\to 0} \frac{h(r)}{g(r)}H^{g}(A).
\end{equation}
Moreover, throughout this paper we denote the gauge functions of form 
\begin{equation}\label{Gaugedef}
h(r)=\left( \frac{1}{\log\left( \frac{1}{r} \right)} \right)^{a},
\end{equation}
where $a>0$, by $h_a$.

\section{Multifractal spectra}\label{MULTI}
\begin{figure}[ht]
\begin{center}
\includegraphics[scale=0.75]{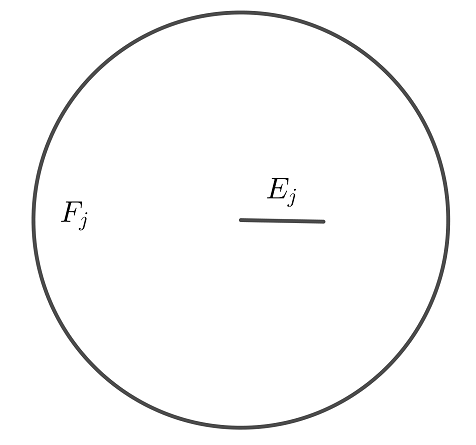}
\caption{}
\label{FIG1}
\end{center}
\end{figure}

The main idea behind the proof of Theorem \ref{MSPEKTRI} is to use capacity estimates which are similar to those in \cite{CH} and \cite{KT}, but have been adapted to the fact that we must measure stretching at many points simultaneously. We will use disjoint sets that consist of the union of  balls $F_j$ and  line-segments $E_j$, see the figure \ref{FIG1}, as  building blocks for our condensers. 

\noindent One of the first obstacles  we encounter is to ensure that we can find sufficiently many disjoint balls $F_j$, such that all of them have approximately the same size and that the line segments $E_j$ inside them satisfy strong stretching properties. To this end, we use the following lemma.
\begin{lemma}\label{KUULAL}
Let $a>0$ and $A \subset \mathbb{R}^n$ be given such that $H^{h_a}(A)>0$, where 
\begin{equation*}
h_a(r)=\left( \frac{1}{\log\left( \frac{1}{r} \right)} \right)^{a}.
\end{equation*}
Furthermore, assume that for every point $z\in A$ there exists a decreasing sequence $r_{z,n} \to 0$. Then for any given $0<b<a$ we can find 
\begin{equation}
\left\lfloor \log^{b}\left( e^{e^{i}} \right) \right\rfloor
\end{equation}
disjoint balls $B\left( z_j, r_{j}^{\frac{1}{2}} \right)$, where $z_j\in A$, $r_j \in \{r_{z_j,n}\}_{n=1}^{\infty}$ and $r_j \in I_i=\left[ \frac{1}{e^{e^{i}}},  \frac{1}{e^{e^{i-1}}} \right)$ for every $j$. Moreover, we can choose the exponent $i$ as big as we wish, and thus the radii $r_j$ can be made arbitrary small.
\end{lemma}
\proof Let us assume that the claim is false, that is, there exists $b<a$ and $i_0$ such that we can not find sufficiently many suitable balls for any $i>i_0$, and derive a contradiction. \newline \newline
Fix an arbitrary $i> i_0$ and denote by $A_i$ the set of points $z\in A$ for which there exists a radius $r_{z} \in \{r_{z,n}\}_{n=1}^{\infty}$ such that $r_z \in \left[ \frac{1}{e^{e^{i}}},  \frac{1}{e^{e^{i-1}}} \right)$. The set $A_i$ might be empty, in which case we move on to the next integer. If the set $A_i$ is non-empty we choose for every point $z\in A_i$ some radius $r_z$ and fix the ball $B\left( z, r_{z}^{\frac{1}{2}} \right)$. Using Vitali's covering theorem we can select countable many disjoint balls $ B\left( z_j, r_{z_j}^{\frac{1}{2}} \right)$ such that
\begin{equation*}
A_i \subset \bigcup_{j\in J_i} B\left( z_j, 5r_{z_j}^{\frac{1}{2}} \right).
\end{equation*}
Moreover, according to our assumption $||J_i|| \leq 2 \log^{b}\left( e^{e^{i}} \right)$ for every $i$. \newline \newline
Since  
\begin{equation*}
A \subset \bigcup_{i>i_0} A_i
\end{equation*}
we can choose a constant $d$ such that $b<d<a$, and estimate the measure $H^{h_d}(A)$ using these balls $B\left( z_j, 5r_{z_j}^{\frac{1}{2}} \right)$. To this end we calculate
\begin{equation*}
\begin{split}
\sum_{\substack{i>i_0, \\  j\in J_i}} h_d\left( \text{diam} \left(B\left( z_j,5r_{z_j}^{\frac{1}{2}} \right) \right) \right) & \leq c_d\sum_{\substack{i>i_0, \\  j\in J_i}} \frac{1}{\log^{d}\left( \frac{1}{r_{z_j}} \right)} \leq c_d\sum_{\substack{i>i_0, \\  j\in J_i}}  \frac{1}{e^{d(i-1)}} \\
& \leq c_{d}\sum_{i\geq 1} \frac{\log^{b}\left( e^{e^{i}} \right)}{e^{d(i-1)}}=  c_{d}\sum_{i\geq 1} \left(\frac{e^b}{e^d}\right)^{i}< \infty,
\end{split}
\end{equation*}
as $b<d$. Thus we see that $H^{h_d}(A)< \infty$. But this is a contradiction with the assumption that $H^{h_a}(A)>0$, due to the inequality \eqref{QEY}, and hence the claim holds.

\subsection{Proof of Theorem \ref{MSPEKTRI}}\label{MULTI.1}
Let us then use Lemma \ref{KUULAL} to prove Theorem \ref{MSPEKTRI}. Here we write the stretching condition \eqref{MSYHT} in the form 
\begin{equation}\label{MSYHT2}
\left|f^{-1}(z+\lambda_{z,n})-f^{-1}(z)\right| \geq c \left( \log\left( \frac{1}{|\lambda_{z,n}|} \right) \right)^{\frac{-p(n-1)-d}{n}},
\end{equation}
where $d>0$ is some given constant, and show that  
\begin{equation*}
H^{h_a}(A)=0
\end{equation*}
for every $a> d$. Theorem \ref{MSPEKTRI} then follows by choosing $d= \frac{ps(n-1)}{n-s}$. Note, that we can additionally assume without loss of generality that $A\subset B(0,1)$.  \newline \newline
Fix $a>d$ and assume for a moment that $H^{h_a}(A)>0$. Our aim is to show that this leads to a contradiction. \newline \newline 
By this assumption Lemma \ref{KUULAL} yields that for any constant $b$, with $d<b<a$, and for some arbitrary big integers $i$ there exists 
\begin{equation*}
\left\lfloor \log^{b}\left( e^{e^{i}} \right) \right\rfloor
\end{equation*}
disjoint balls $B_j=B\left( z_j, r_{z_j}^{\frac{1}{2}} \right)$, where $z_j \in A$, $r_{z_j}\in \{ |\lambda_{z_j,n}| \}_{n=1}^{\infty}$ and $r_{z_j} \in I_i= \left[ \frac{1}{e^{e^{i}}},  \frac{1}{e^{e^{i-1}}} \right)$ for every $j$. Let us denote the union of these balls by $F$ and the union of the line segments $E_j=\left[ z_j, z_j+ \lambda_{z_j,n_j} \right]$, where $r_{z_j}=|\lambda_{z_j,n_j}|$, by $E$. \newline \newline
The pairs $E$, $F$ and $f^{-1}(E)$, $f^{-1}(F)$ form condensers, see the figure \ref{FIG1}, and estimates for their capacities will play a central role in the proof. As the balls $B_j$ are disjoint and the mapping $f$ is a homeomorphism these capacities can be calculated as the sum of the capacities for the condensers formed by the  pairs $E_j, B_j $ and $f^{-1}(E_j), f^{-1}(B_j)$.  \newline \newline
We start by fixing  $q=\frac{pn}{p+1}$ and estimate the capacity 
\begin{equation*}
\text{cap}_{q}\left(f^{-1}(E),f^{-1}(F)\right)
\end{equation*}
 from below. Note that since $p>n-1$ it holds that $n-1<q<n$. Thus we can use Lemma \ref{KAPLEM} and the stretching estimate \eqref{MSYHT2} to obtain
\begin{equation}\label{APUYHT1}
\begin{split}
\text{cap}_{q}\left( f^{-1}(E),f^{-1}(F) \right)  &=  \sum_{j} \text{cap}_{q}\left( f^{-1}(E_j), f^{-1}(B_j) \right) \\
& \geq   c_{n,p} \sum_{j} |f^{-1}(z_j)-f^{-1}(z_j+\lambda_{z_j,n_j})|^{n-q} \\
& \geq c_{n,p} \sum_{j} \left( \log \left( \frac{1}{r_{z_j}} \right) \right)^{\frac{-p(n-1)-d}{p+1}} \\
& \geq  c_{n,p} \sum_{j} \left( \log \left( e^{e^{i}} \right) \right)^{\frac{-p(n-1)-d}{p+1}} \\
& \geq c_{n,p} \left( \log \left( e^{e^{i}} \right) \right)^{b-\frac{p(n-1)+d}{p+1}}.
\end{split}
\end{equation}
Next we use the observation \eqref{KAPTULOS} to estimate the capasity
\begin{equation*}
\text{cap}_{n}(E,F)
\end{equation*}
from above by
\begin{equation}\label{APUYHT2}
\begin{split}
\text{cap}_{n}(E,F)  = \sum_{j}\text{cap}_{n}(E_j,F_j) & \leq c_n \sum_{j}  \log^{1-n}\left( \frac{r_{j}^{\frac{1}{2}}}{r_j} \right) \\
& \leq c_n \sum_{j} \left(\log\left( e^{\frac{1}{2}e^{i-1}} \right)\right)^{1-n} \\
& = c_n \sum_{j} \left(\log\left( e^{\frac{1}{2e}e^{i}} \right)\right)^{1-n} \\
& \leq c_n \sum_{j} \left( \log \left( e^{e^{i}} \right) \right)^{1-n} \\
& \leq c_{n}  \left( \log \left( e^{e^{i}} \right) \right)^{b+1-n}.
\end{split}
\end{equation}
Finally, we provide a relation between these capacities, in the spirit of \cite{KT} and \cite{CH}, and show that the stretching condition \eqref{MSYHT2} can only be satisfied in a small set. \newline \newline
Let $u$ be an admissible function for the condenser $(E,F)$. Set $v=u \circ f$, and note that since $f$ is a homeomorphism $v$ is admissible for the condenser $\left(  f^{-1}(E), f^{-1}(F) \right)$. From the chain rule and the distortion inequality we obtain 
\begin{equation*}
|\nabla v(z)|^{q}\leq|\nabla u(f(z))|^{q}K_{f}^{\frac{q}{n}}(z)J_{f}^{\frac{q}{n}}(z).
\end{equation*}
Hence we can use H\"older's inequality and a change of variables to estimate
\begin{equation}\label{APUYHT3}
\begin{split}
\int_{f^{-1}(F)} |\nabla v(z)|^{q} \: dz & \leq \int_{f^{-1}(F)} |\nabla u(f(z))|^{q}K_{f}^{\frac{q}{n}}(z)J_{f}^{\frac{q}{n}}(z) \: dz \\
& \leq \left( \int_{f^{-1}(F)} \left(K_{f}^{\frac{q}{n}}(z)\right)^{\frac{n}{n-q}} \; dz \right)^{\frac{n-q}{n}} \left( \int_{f^{-1}(F)}  |\nabla u(f(z))|^{n}J_{f}(z) \; dz\right)^{\frac{q}{n}} \\
& \leq  \left( \int_{f^{-1}(F)} K_{f}^{p}(z) \; dz\right)^{\frac{1}{p+1}} \left( \int_{F} |\nabla u(y)|^{n} \; dy \right)^{\frac{p}{p+1}} \\
& = c_{p,f}  \left( \int_{F} |\nabla u(y)|^{n } \; dy \right)^{\frac{p}{p+1}},
\end{split}
\end{equation}
where in the last equality we have used the fact that $f^{-1}(F) $ lies inside some  compact set. Taking infimum over all admissable functions $u$ we obtain
\begin{equation*}
\text{cap}_{q}^{p+1}\left( f^{-1}(E),f^{-1}(F) \right) \leq c_{p,f} \cdot \text{cap}_{n}^{p}(E,F).
\end{equation*}
Combining this with the estimates \eqref{APUYHT1} and \eqref{APUYHT2} for the capacities we obtain
\begin{equation*}
c_{n,p} \left( \log \left( e^{e^{i}} \right) \right)^{b(p+1)-p(n-1)-d} \leq c_{n,p,f}  \left( \log \left( e^{e^{i}} \right) \right)^{pb+p-np},
\end{equation*}
which simplifies to
\begin{equation*}
\left( \log \left( e^{e^{i}} \right) \right)^{b-d} \leq c_{n,p,f}.
\end{equation*}
But since $b>d$ this can not hold for big $i$, and hence we arrive at a contradiction. Thus the assumption that $H^{h_a}(A)>0$ for some $a>d$ is false and Theorem \ref{MSPEKTRI} holds.

\section{Area compression}\label{AREA}
After establishing Theorem \ref{MSPEKTRI} we can turn our attention to area compression. In order to utilize the stretching multifractal spectra in the  proof of Theorem \ref{AREAC} we need the following lemma, that lets us to partition the general case into suitable pieces.
\begin{lemma}\label{Siivu}
Fix $s\in (0,n)$, $\epsilon>0$, $b>0$ and $p> n-1$, and let $f:\mathbb{R}^n \to \mathbb{R}^n$ be a homeomorphism with $p$-integrable distortion. Assume furthermore, that for every point $z \in A \subset \mathbb{R}^n$ there exists a sequence $\lambda_{z,n}$, such that $|\lambda_{z,n}|\to 0$, for which
\begin{equation}\label{Siivu1}
|f^{-1}(z+\lambda_{z,n})-f^{-1}(z)| \geq \left( \frac{1}{\log \left( \frac{1}{|\lambda_{z,n}|}  \right)} \right)^{\frac{p(n-1)}{n-s}},
\end{equation}
but that all sufficiently small $|\lambda| < r_z$ satisfy 
\begin{equation}\label{Siivu2}
|f^{-1}(z+\lambda)-f^{-1}(z)| \leq b\left( \frac{1}{\log \left( \frac{1}{|\lambda|}  \right)} \right)^{\frac{p(n-1)}{n-(s-\epsilon)}}.
\end{equation}
Then $$\dim(f^{-1}(A))\leq s+g(\epsilon),$$ where $g(\epsilon) \to 0$ when $\epsilon \to 0$.
\end{lemma}
\proof  By Theorem \ref{MSPEKTRI} we know that 
\begin{equation*}
H^{h_a}(A)=0, \quad \text{when} \quad a=\frac{p(s+\epsilon)(n-1)}{n-s-\epsilon} > \frac{ps(n-1)}{n-s}.
\end{equation*}
Since $H^{h_a}(A)=0$ we can find balls $B(z_j,r_j)$ such that $z_j\in A$ for every $j$, the diameter $r_j$ is small enough so that the condition \eqref{Siivu2} holds inside the ball, the union of the balls satisfies 
\begin{equation*}
A \subset \bigcup_{j} B(z_j,r_j)
\end{equation*}
and finally that
\begin{equation*}
\sum_{j} h_a\left( \text{diam}\left( B(z_j,r_j) \right) \right)< \bar{\epsilon},
\end{equation*}
where the constant $\bar{\epsilon}>0$ can be chosen as small as we wish. \newline \newline
Using the images $f^{-1}\left( B(z_j,r_j) \right)$ of these balls with the stretching bound \eqref{Siivu2} we can estimate the 
\begin{equation*}
(s+\epsilon)\cdot \frac{n-s+\epsilon}{n-s-\epsilon}- \text{dimensional}
\end{equation*}
Hausdorff measure of the set $f^{-1}(A)$ by  
\begin{equation*}
\begin{split}
\sum_{j} \left(\text{diam}\left( f^{-1}\left( B(z_j,r_j) \right) \right)\right)^{(s+\epsilon)\cdot \frac{n-s+\epsilon}{n-s-\epsilon}} & \leq c_{b,s,n,\epsilon} \sum_{j} \left( \frac{1}{\log\left( \frac{1}{r_j} \right)} \right)^{\frac{p(s+\epsilon)(n-1)}{n-s-\epsilon}} \\
& \leq c_{b,s,n,\epsilon} \sum_{j} h_a\left( \text{diam}\left( B(z_j,r_j) \right) \right) \\
& <c_{b,s,n,\epsilon} \cdot \bar{\epsilon}.
\end{split}
\end{equation*}
This shows that $\dim\left( f^{-1}(A) \right) \leq (s+\epsilon)\cdot \frac{n-s+\epsilon}{n-s-\epsilon}$, and it is easy to see that this dimension has the correct form of $s+g(\epsilon)$, where $g(\epsilon) \to 0$ as $\epsilon \to 0$. \newline

\subsection{Proof of Theorem \ref{AREAC}}\label{AREACS}
With  Lemma \ref{Siivu} at our disposal we can proceed to prove Theorem \ref{AREAC} on  area compression. \newline \newline
To this end, we will show that if we fix $s\in (0,n)$ and assume that a set $A \subset D$ satisfies
\begin{equation*}
H^{h_a}(A)=0, \quad \text{where} \quad a=\frac{ps(n-1)}{n-s},
\end{equation*}
then 
\begin{equation*}
H^{\bar{s}}\left(f^{-1}(A)\right)=0
\end{equation*}
for every $\bar{s}>s$. \newline \newline
So, let us fix some $\bar{s}>s$ and let $A_0 \subset A$ be the set of those points $z\in A$ for which there exists radius $r_z$ such that
\begin{equation*}
|f^{-1}(z+\lambda)-f^{-1}(z)|\leq \left( \frac{1}{\log \left( \frac{1}{|\lambda|} \right)}\right)^{\frac{p(n-1)}{n-s}},
\end{equation*}
when $|\lambda|< r_{z}$. Then we can use the fact that $H^{h_a}(A_0)=0$ to choose balls $B(z_j,r_j)$ in a similar manner as in the previous lemma and estimate
\begin{equation*}
\sum_{j} \text{diam}^{s}\left( f^{-1}(B(z_j,r_j)) \right) \leq c_s \sum_{j} \left( \frac{1}{\log \left( \frac{1}{r_j} \right)}\right)^{\frac{ps(n-1)}{n-s}}<c_s\bar{\epsilon}.
\end{equation*}
Particularly, we see that $H^{\bar{s}}(f^{-1}(A_0))=0$. \newline \newline
Then we start using Lemma \ref{Siivu}. First, choose $\epsilon$ such that
\begin{equation*}
(s+\epsilon)\cdot \frac{n-s+\epsilon}{n-s-\epsilon}<\bar{s}
\end{equation*}
and fix $s_1=s-\epsilon$. Then, denote by $A_1$ the set of points $z\in A$  for which there exists  a sequence $\lambda_{z,n}$, satisfying $|\lambda_{z,n}| \to 0$ when $n\to \infty$, such that
\begin{equation*}
|f^{-1}(z+\lambda_{z,n})-f^{-1}(z)| \geq \left( \frac{1}{\log \left( \frac{1}{|\lambda_{z,n}|}  \right)} \right)^{\frac{p(n-1)}{n-s}},
\end{equation*}
but
\begin{equation*}
|f^{-1}(z+\lambda)-f^{-1}(z)| \leq \left( \frac{1}{\log \left( \frac{1}{|\lambda|}  \right)} \right)^{\frac{p(n-1)}{n-s_1}}
\end{equation*}
for all sufficiently small $|\lambda|<r_z$. Then Lemma \ref{Siivu} asserts that
\begin{equation*}
H^{(s+\epsilon)\cdot \frac{n-s+\epsilon}{n-s-\epsilon}}(f^{-1}(A_1))<\infty,
\end{equation*}
and thus  
\begin{equation*}
H^{\bar{s}}(f^{-1}(A_1))=0.
\end{equation*}
Next choose $s_2=s_1-\epsilon$ and note that 
\begin{equation}\label{aaaa}
 (s_1+\epsilon)\cdot \frac{n-s_1+\epsilon}{n-s_1-\epsilon}<(s+\epsilon)\cdot \frac{n-s+\epsilon}{n-s-\epsilon}<\bar{s}.
\end{equation}
Then we use Lemma \ref{Siivu} again  and define the set $A_2$ to be the set of points $z\in A$  for which there exists  a sequence $\lambda_{z,n}$, satisfying $|\lambda_{z,n}| \to 0$ when $n\to \infty$, such that
\begin{equation*}
|f^{-1}(z+\lambda_{z,n})-f^{-1}(z)| \geq \left( \frac{1}{\log \left( \frac{1}{|\lambda_{z,n}|}  \right)} \right)^{\frac{p(n-1)}{n-s_1}},
\end{equation*}
but
\begin{equation*}
|f^{-1}(z+\lambda)-f^{-1}(z)| \leq \left( \frac{1}{\log \left( \frac{1}{|\lambda|}  \right)} \right)^{\frac{p(n-1)}{n-s_2}}
\end{equation*}
for all sufficiently small $|\lambda|<r_z$. Then Lemma \ref{Siivu} with the inequality \eqref{aaaa} implies 
\begin{equation*}
H^{\bar{s}}(f^{-1}(A_2))=0.
\end{equation*}
We continue in a similar manner, using the fact that
\begin{equation*}
 (s_{i+1}+\epsilon)\cdot \frac{n-s_{i+1}+\epsilon}{n-s_{i+1}-\epsilon}<(s_i+\epsilon)\cdot \frac{n-s_i+\epsilon}{n-s_i-\epsilon}<\bar{s}
\end{equation*}
at every step $i$, until we choose $s_n=s_{n-1}-\epsilon=0$. We can guarantee that such $s_n$ exists by choosing suitable $\epsilon$. \newline \newline 
Finally, we define the set $A_n$ to consist of points $z\in A$  for which there exists  a sequence $\lambda_{z,n}$, satisfying $|\lambda_{z,n}| \to 0$ when $n\to \infty$, such that
\begin{equation*}
|f^{-1}(z+\lambda_{z,n})-f^{-1}(z)| \geq \left( \frac{1}{\log \left( \frac{1}{|\lambda_{z,n}|}  \right)} \right)^{\frac{p(n-1)}{n-s_{n-1}}},
\end{equation*}
but
\begin{equation*}
|f^{-1}(z+\lambda_{z})-f^{-1}(z)| \leq c_{f,p,n} \left( \frac{1}{\log \left( \frac{1}{|\lambda_{z}|}  \right)} \right)^{\frac{p(n-1)}{n}}
\end{equation*}
for all sufficiently small $|\lambda|<r_z$. Then Lemma \ref{Siivu} with the inequality \eqref{aaaa} implies 
\begin{equation*}
H^{\bar{s}}(f^{-1}(A_n))=0.
\end{equation*}
The modulus of continuity result \eqref{CHyht}, where we choose $G=f^{-1}(B(0,2))$ for every point $z\in A$, verifies that 
\begin{equation*}
A=\bigcup_{i=0}^{n}A_i,
\end{equation*}
and thus 
\begin{equation*}
H^{\bar{s}}(f^{-1}(A)) \leq H^{\bar{s}} \left( \bigcup_{i=0}^{n} f^{-1}(A_i) \right)=0.
\end{equation*}
This finishes the proof of Theorem \ref{AREAC}.


\begin{thebibliography}{10}
\bibitem{AS} K. Astala, \textit{Area distortion of quasiconformal mappings,} Acta Math., 173 (1994), 37-60. 

\bibitem{ACM} K. Astala, A. Clop, J. Mateu, J. Orobitg, and I. Uriarte-Tuero, \textit{Distortion of Hausdorff measures and improved Painlev\'{e} removability for bounded
quasiregular mappings,} Duke Math. J., 141 (2008), 539-571.

\bibitem{AIM}  K. Astala, T. Iwaniec, and G. J. Martin, \textit{Elliptic partial differential equations and quasiconformal mappings in the plane}, Princeton University Press, 2009.



\bibitem{CH}  A. Clop and D. Herron, \textit{Mappings with finite distortion in $L^{p}_{\text{loc}}$: Modulus of continuity and compression of Hausdorff measure}, D.A. Isr. J. Math. (2014) 200: 225. 


\bibitem{HKM} J. Heinonen, T. Kilpel\"{a}inen, and O. Martio, \textit{Nonlinear potential theory of degenerate elliptic equations}, Oxford Univ. Press, Oxford, 1993.

\bibitem{HK} S. Hencl and P. Koskela, \textit{Lectures on mappings of finite distortion}, Lecture Notes in Mathematics, vol. 2096, Springer, Cham, 2014.


\bibitem{H4} L. Hitruhin, \textit{Joint rotational and stretching multifractal spectra of mappings with integrable distortion}, to appear in Revista Matem\'{a}tica Iberoamericana.




\bibitem{KT} P. Koskela and J. Takkinen, \textit{Mappings of finite distortion: formation of cusps. III.} Acta Math. Sin. (Engl. Ser.) 26(5), 817-824 (2010).




\bibitem{U} I. Uriarte-Tuero,  \textit{Sharp examples for planar quasiconformal distortion of Hausdorff measures and removability}, Int. Math. Res. Notices IMRN, (2008), 43 pp.


\bibitem{VU} M. Vuorinen, \textit{Conformal geometry and quasiregular mappings,} Lecture Notes in Math., 1319, Springer-Verlag, Berlin-New York, 1988. 





\end{thebibliography}
\end{document}